\documentclass[final,5p,times,twocolumn]{elsarticle}

\usepackage{graphicx}
\usepackage{subcaption}
\usepackage{amssymb, amsmath}
\usepackage{amsthm}
\usepackage{color}
\usepackage{array}
\usepackage{tabu}
\usepackage{enumitem}

\usepackage[utf8]{inputenc}
\usepackage[english]{babel}
\usepackage[T1, T2A]{fontenc}

\biboptions{compress}

\newtheorem{theorem}{Theorem}[section]

\theoremstyle{definition}
\newtheorem{definition}[theorem]{Definition}

\theoremstyle{remark}
\newtheorem{remark}[theorem]{Remark}

\begin{document}

\begin{frontmatter}

\title{On the Application of the Wa{\.z}ewski Method\\ to the Problem of  Global Stabilization}



\author[label5]{Ivan Polekhin}
\address[label5]{Steklov Mathematical Institute of the Russian Academy of Sciences, Moscow, Russia}


\begin{abstract}
We consider a possible application of the Ważewski topological method to feedback control systems and to more general dynamical systems. We show how this method can be used to prove the impossibility of global stabilization in such problems. Moreover, we give sufficient conditions for the existence of a solution such that its trajectory never leaves a subset of the extended phase space of the system and does not tend asymptotically to a given equilibrium. We illustrate our result with various real-life systems including the Furuta pendulum and the wheeled inverted pendulum.
\end{abstract}


\end{frontmatter}

\section{Introduction} 		      

The design of feedback control is among the major problems of applied mathematics. In applications we often face the task of making some configuration of the system asymptotically stable in the sense of Lyapunov. For this we can usually use a feedback control, and this control is constrained by the general design of our system. At the same time, we are not only interested in stability, but also it is often required to make the corresponding basin of attraction as large as possible.

In particular, one can try to find sufficient conditions under which the system is stable in the large and all its solutions tend to a unique equilibrium. For equations defined on $\mathbb{R}^n$, an overview of such results can be found in \cite{krasovskii1963stability} (see also \cite{hartman1961stability,garg1989global,lohmiller1998contraction})

For systems defined on a closed manifold, a fundamental result was proved in \cite{bhat2000topological}. To be more precise, it was shown that, given a continuous semi-flow on a manifold $M$ such that there exists a vector bundle $\pi \colon M \to N$, with $N$ a closed manifold, it is impossible for the system to have a globally asymptotically stable equilibrium. In particular, if we consider a dynamical system acting on the tangent bundle $TN$ of a closed manifold $N$, then this system cannot have such an equilibrium.

As an example of such a system one can consider a planar pendulum with feedback control: the system defined on $T\mathbb{S}^1$ cannot be globally asymptotically stable. However, if we impose constraints on this pendulum system, the phase space of the corresponding system can change. If the pendulum is placed on a horizontal plane of support, then its phase space is $\mathbb{R} \times [0, \pi]$ and the result from \cite{bhat2000topological} can't be applied here. It is also possible to consider systems defined not by semi-flows, but by so-called semi-processes. For instance, this will be the case when the corresponding system of ODEs is non-autonomous, i.e., explicitly depends on time.

The Ważewski topological method \cite{wazewski1947principe,reissig1963qualitative,srzednicki2004wazewski} has been already considered by many authors as a useful tool in nonlinear analysis. It can be applied to general problems when one needs to determine the asymptotic behavior of the solutions of a given system \cite{reissig1963qualitative,conley1975application,onuchic1961applications} as well as to more specific problems, including finding periodic solutions \cite{srzednicki1994periodic,srzednicki2005fixed}, solving boundary value problems \cite{conley1975application} or even detecting chaotic behavior of the system \cite{wojcik2011wazewski}.

The main idea of the Ważewski topological method can be also applied if we want to prove the impossibility of global stabilization in systems with feedback control. To be more precise, we will show that if a uniformly Lyapunov stable equilibrium is located in a larger set such that solutions of the system intersect its boundary transversely (or, in a more general case, the set of strict egress points coincides with the set of egress points, if we use the terminology of the Ważewski method), then there is a solution with the following properties:
\begin{itemize}
    \item This solution does not tend to an equilibrium as time tends to infinity.
    \item Its trajectory never leaves the considered set.
\end{itemize}

In this paper, we present the above more formally. The present paper can be considered as a development of the result in \cite{polekhin2018topological}: we remove the requirement that there exists a Lyapunov function in a vicinity of the equilibrium and present a general result for a semi-process, not limiting ourselves to the study of a specific system. We also present some mechanical examples. Comparing with the result in \cite{bhat2000topological}, we might say that \cite{bhat2000topological} gives a solution to the problem of global stabilization on a manifold without boundary and in our paper we outline a possible approach to the same problem on a manifold with boundary.

This paper is organized as follows. First, we briefly explain the main theorem of the Ważewski topological method for flows. Then we present our main result on the impossibility of global stabilization. The proof is self-contained and can be understood without any external references. The result is illustrated by three classical control systems: the inverted pendulum, the wheeled inverted pendulum, and the Furuta pendulum.

\section{Results}
\subsection{The Wa{\.z}ewski method for flows}

We start with a brief explanation of the Wa{\.z}ewski topological method. First, we consider the case of flows on a smooth manifold.  Let $M$ be a smooth ($C^\infty$) manifold. Let $\varphi \colon M \times \mathbb{R} \to M$ be a flow on $M$, i.e., it is a continuous map such that 
\begin{enumerate}
    \item For any $x \in M \times \mathbb{R}$, we have $\varphi(x,0) = x$,
    \item For any $t, s \in \mathbb{R}$, we have $$\varphi(x,t+s) = \varphi(\varphi(x,t),s).$$
\end{enumerate}
Let $W \subset M$ be an open subset of $M$ with non-empty boundary: $\partial W \ne \varnothing$. Note that $W$ can be an arbitrary set, but in the applications it is usually a set with a piecewise smooth boundary.

For the Wa{\.z}ewski method, the following two notions of egress and strict egress points play the key r{\^o}le. 
\begin{definition}
We say that point $x \in \partial W$ is an egress point for $W$ w.r.t. the flow $\varphi$ if there exists an $\varepsilon > 0$ such that $\varphi(x,t) \in W$ for all $t \in (-\varepsilon,0)$.
\end{definition}
\begin{definition}
We say that the egress point $x \in \partial W$ is a strict egress point for $W$ w.r.t. the flow $\varphi$ if there exists an $\varepsilon > 0$ such that $\varphi(x,t) \notin W \cup \partial W$ for all $t \in (0,\varepsilon)$.
\end{definition}
The set of all egress points will be denoted by $W^-$ and the set of all strict egress points is $W^{--}$. Roughly speaking, the condition $W^- = W^{--}$, which we are going to use below, means that there are no trajectories internally tangent to $\partial W$, i.e., all solutions are either transverse to the boundary or externally tangent to it.

For the case of smooth ODEs, one can check whether a point is an egress point or a strict egress point by examining the corresponding Taylor expansion at the point.

For any point $x \in M$, let us consider the half-trajectory of the flow:
$$
\gamma_\tau(x) = \bigcup_{t \in [0,\tau)} \varphi(x,t) \subset M.
$$
\begin{definition}
For $x \in W$, we say that $$\sigma(x) = \sup \{ \tau \geqslant 0 \colon \gamma_\tau(x) \subset W \}.$$  is the time of egress from $W$. For $x \in \partial W$ we put $\sigma(x) = 0$. If $\sigma(x) = \infty$, we say that the half-trajectory starting at $x$ does not leave $W$.
\end{definition}
The Wa{\.z}ewski method provides a robust approach to proving that there exist points $x \in W$ such that $\sigma(x) = \infty$, i.e. the corresponding trajectories never leave $W$.
\begin{theorem}
Let $M$ be a manifold and $\varphi \colon M \times \mathbb{R} \to M$ be a flow. If there exists an open set $W \subset M$ and a set $\Gamma \subset W \cup W^{-}$ such that $W^{-} = W^{--}$, and $\Gamma \cap W^{-}$ is a retract of $W^{-}$ but is not a retract of $\Gamma$. Then there is a point $x \in \Gamma$ such that $\sigma(x) = \infty$.
\end{theorem}

A proof of this result can be found, for instance, in \cite{hartman1982ordinary} (see also \cite{srzednicki2005fixed,mohamed2016methode}). However, the general idea of the proof can be explained without going into technical details. First, let us recall the definition of a retract. Let $X$ be a topological space and $Y \subset X$. Then $Y$ is a retract of $X$ if there exists a continuous map $r \colon X \to Y$ such that $r(x) = x$ for any $x \in Y$. Equivalently, $Y$ is a retract of $X$ if and only if every continuous mapping of $Y$ into an arbitrary topological space $Z$ can be extended to a continuous mapping of the entire space $X$ into $Z$.

Suppose that for any $x \in \Gamma$, we have $\sigma(x) < \infty$, i.e., all solutions starting in the considered set reach the boundary $\partial W$. The key observation that we are going to use is that $\sigma(x) \colon \Gamma \to \mathbb{R}$ is a continuous function. This follows from the assumption that $W^- = W^{--}$ and the continuity of the flow. Hence, the map $x \mapsto \varphi(x, \sigma(x))$ is also continuous. We will denote this map by $m$ and the retraction between $W^{-}$ and $\Gamma \cap W^{-}$ by $r$. Then $r \circ m$ defines a retraction between $\Gamma$ and $\Gamma \cap W^{-}$ ($m$ maps $\Gamma$ into $W^-$). This contradiction proves the statement. Note that the existence of the retraction $r$ is an assumption of the theorem. That is, we assume that the required retraction exists, and when one applies Theorem 2.4 to some system, the existence of such a retraction has to be proved.

\begin{figure}[h!]
  \centering
    \includegraphics[width=0.47\textwidth]{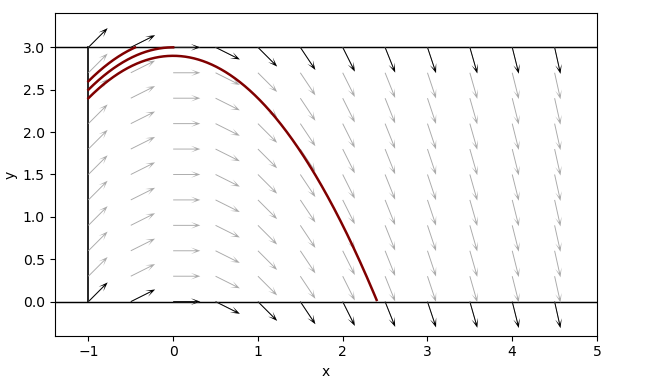}
\caption{All solutions leave $W$.}
\label{fig1}
\end{figure}

Note that the condition $W^- = W^{--}$ is important if we want to construct a continuous map from $\Gamma$ into $W^{-}$. Consider a toy system
\begin{align*}
    & \dot x = 1, \quad \dot y = -x.
\end{align*}
For this system and $W = \{ x,y \colon 0 < y < 3 \}$, we have $W^- \ne W^{--}$ (see Fig. \ref{fig1}), i.e., we have a trajectory that is internally tangent to $\partial W$, therefore, $(0,3)$ of the boundary is an egress point, yet not a strict egress point. Along this trajectory the map $m$ is discontinuous: points located slightly above it leave $W$ through the upper part of the boundary; in contrast, points that are below this curve leave $W$ through the lower part of the boundary. 

\begin{figure}[h!]
  \centering
    \includegraphics[width=0.45\textwidth]{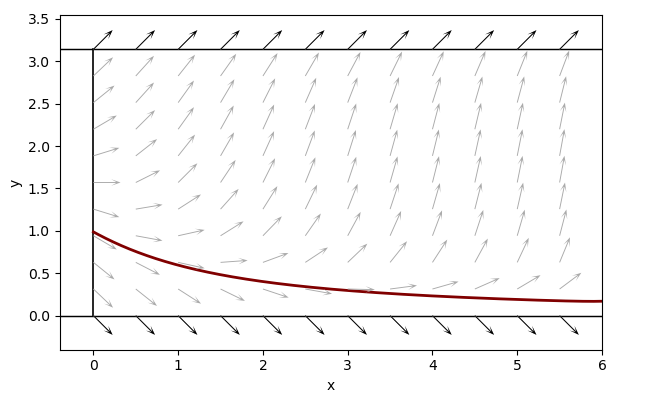}
\caption{A numerically obtained solution that never leaves $W$.}
\label{fig2}
\end{figure}

Before moving on to our main results, we consider a couple of examples illustrating the method. We start with an archetypical example: a vector field that transversely intersects  an infinite strip. For instance, let us consider the nonlinear system
\begin{align*}
    & \dot x = 1,\\
    & \dot y = a \cdot \cos y + x\sin y.
\end{align*}
For $a = 1$, there are infinitely many solutions $(x(t), y(y))$ such that $y(t) \in [0, \pi]$ for all $t \geqslant 0$. This follows from the two inequalities $\dot y \big|_{y = 0} > 0 $ and $\dot y \big|_{y = \pi} < 0$. For $a = -1$ we have $\dot y \big|_{y = 0} < 0 $ and $\dot y \big|_{y = \pi} > 0$, i.e., all solutions of the system leave the strip $W$ defined by the inequalities $0 \leqslant y \leqslant \pi$ transversely to its boundary (see Fig. \ref{fig2}). Therefore, $W^- = W^{--}$ and we can apply Theorem 2.4. Consider the segment $\Gamma$ defined by $x = 0$ and $0 \leqslant y \leqslant \pi$. The set $\Gamma \cap W^-$ consists of two points, and the two lines $W^-$ can be retracted to these points (this retraction is unique). At the same time, $\Gamma$ cannot be retracted to its boundary (a line segement cannot be retracted into its endpoints). Therefore, there is a solution that never leaves $W$. Again, note that the existence of the retraction is an important assumption in the above considerations.

\begin{figure}[h!]
  \centering
    \includegraphics[width=0.45\textwidth]{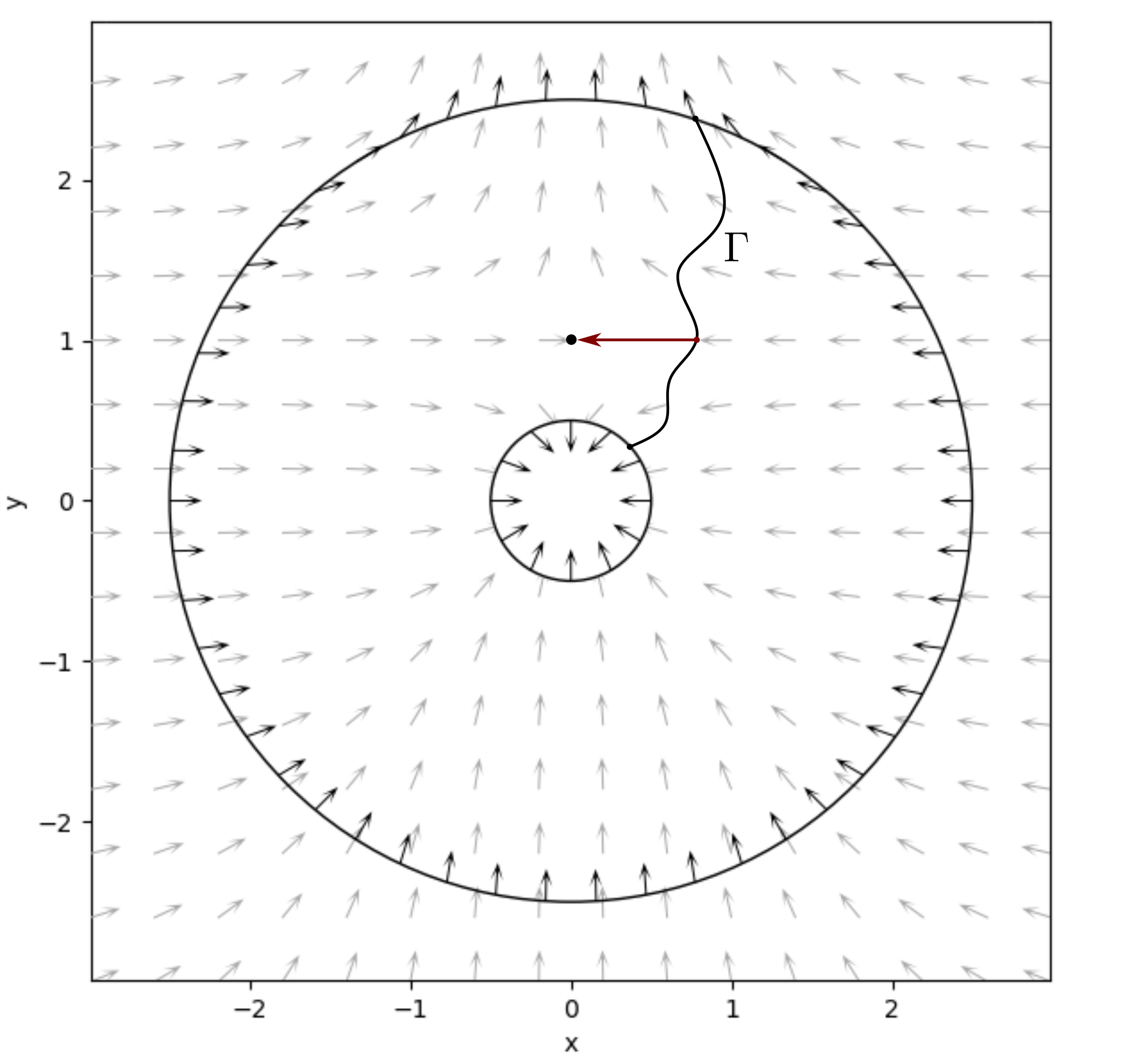}
\caption{The hyperbolic equilibrium $(0,1)$ and the corresponding asymptotic solution are highlighted.}
\label{fig3}
\end{figure}

Another example is the system
\begin{align}
\begin{split}
\label{exam3}
    & \dot x = -x - x^3,\\
    & \dot y = - y + y^2.
\end{split}
\end{align}
Now consider the set $W$ defined by
$$
W = \{ x, y \colon \rho_1 < x^2 + y^2 < \rho_2 \}.
$$
Here, $\rho_1$ is a relatively small parameter and, conversely, $\rho_2$ is a large number (Fig. \ref{fig3}). Then $W$ satisfies the conditions of Theorem 2.4. The part of $\partial W$ defined by $x^2 + y^2 = \rho_1$ consists of strict egress points only. Let $\Gamma$ be an arbitrary smooth curve connecting two disjoint components of $W^{--}$. Then there is a point $x \in \Gamma$ such that the trajectory starting at this point never leaves $W$.

Note that for this particular system, the corresponding solution can be found explicitly: any smooth curve connecting the small circle with the set of egress points of the larger circle inevitably intersects the line $y = 1$, the stable manifold for the hyperbolic equilibrium $x=0$, $y=1$. However, we only used the information about the vector field in a vicinity of the boundary $\partial W$ to prove the existence of such a solution and our approach can be carried over to higher dimensional and non-autonomous cases, where the structure of equilibrium points and invariant manifolds is, in general, unknown.

\subsection{Main results}

Consideration of a system defined by a flow has its advantages, yet this requirement is too restrictive and below we will consider a more general setting that allows us to deal with the following two cases:
\begin{itemize}
    \item Our system can be non-autonomous. In other words, when we consider a system of ordinary differential equations that describes the motion of our feedback control system, we allow an explicit dependence on time for the right hand side. For instance, this dependence can be considered as an external (uncontrollable and unavoidable) force acting on our system.
    \item The right hand side of our ODE can be less regular than is required for the uniqueness of solutions. Actually, what we will need is the so-called right-uniqueness of the solutions.
\end{itemize}

Both these generalizations can be embraced by the notion of a semi-process. Let $\Phi$ be a continuous semi-process (see, for instance, \cite{srzednicki2005fixed}):
\begin{align}
\label{eq1}
\Phi \colon M \times \mathbb{R} \times [0, \infty) \to M.
\end{align}
In other words, $\Phi$ is a continuous map such that
$$
\varphi \colon (x, t_0,t) \mapsto (\Phi(x,t_0,t), t_0 + t) \in M \times \mathbb{R}
$$
is a continuous semi-flow on $M \times \mathbb{R}$, i.e.:
\begin{enumerate}
    \item A map $\varphi \colon M \times \mathbb{R} \times [0, \infty) \to M \times \mathbb{R}$ which is continuous;
    \item For any $(x, t_0) \in M \times \mathbb{R}$, we have $\varphi(x,t_0,0) = (x,t_0)$;
    \item For any $t, s \in [0, \infty)$, we have $$\varphi(x,t_0,t+s) = \varphi(\varphi(x,t_0,t),s).$$
\end{enumerate}

Below we will use the notation $$\Phi(x,t_0,t) = \Phi_{t_0, t}(x).$$
Note that any continuous semi-flow on $M$ can be considered as a continuous semi-process on $M \times \mathbb{R}$ with no dependence on $t_0$.

We now introduce the notions of egress and strict egress points for semi-processes. Note that semi-flows and semi-processes are defined only for non-negative values of $t$. Therefore, it is impossible to carry over Definition 2.1 directly to the case of a semi-process.

First, similarly to the case of a flow, we can consider the half-trajectory of the semi-process:
\begin{align}
    \gamma_\tau(x) = \bigcup_{t \in [0,\tau)} (\Phi_{0,t}(x),t) \subset M \times \mathbb{R}.
\end{align}
From now on, we fix the initial moment of time to be zero. 

Everywhere below we assume that $W \subset M \times \mathbb{R}$ is such a set that $W \cap \{ t = 0 \} \ne \varnothing$, i.e., set $W$ has a non-empty intersection with the plane $t = 0$.

\begin{definition}
Given a point $(x,0) \in W$, we say that $$\sigma(x) = \sup \{ \tau \geqslant 0 \colon \gamma_\tau(x) \subset W \}.$$  is the time of egress from $W$. For $(x,0) \in \partial W$ we put $\sigma(x) = 0$.
\end{definition}
\begin{definition}

We say that $(x_1, t_1) \in \partial W$, $t_1 > 0$ is an egress point for $W$ if there exists a point $(x,0) \in W$ such that
\begin{align}
    (x_1, t_1) = (\Phi_{0,\sigma(x)}(x),\sigma(x)).
\end{align}
Since we cannot consider our system in reverse time, but can only consider our system for $t \geqslant 0$, we have to deal with the case $t_1 = 0$ separately. A point $(x_1, t_1) \in \partial W$, $t_1 = 0$ is an egress point if for some $\varepsilon > 0$ we have $$(\Phi_{0,t}(x_1),t) \notin W \cup \partial W$$ for all $t \in (0, \varepsilon)$.
We denote the set of all egress points by $W^-$.
\end{definition}

\begin{definition}
 We say that an egress point $(x_1, t_1) \in \partial W$ is a strict egress point for $W$ if for some $\varepsilon > 0$ we have $$(\Phi_{0,\sigma(x)+t}(x),\sigma(x)+t) \notin W \cup \partial W$$ for all $t \in (0, \varepsilon)$. Here $(x_1, t_1) = (\Phi_{0,\sigma(x)}(x),\sigma(x))$. We denote the set of all strict egress points by $W^{--}$.
\end{definition}

\begin{remark}
We call a point $(x_1, t_1) \in \partial W$, $t_1 = 0$ an egress point even when this point is actually a strict egress point. We do this only for technical reasons. This will not cause any ambiguity, since everywhere below we assume that for our systems $W^- = W^{--}$.
\end{remark}


\begin{definition}
We say that $x_0$ is an equilibrium for the semi-process (\ref{eq1}) if $\Phi_{t_0, t}(x_0)=x_0$ for all $t_0 \geqslant 0$ and  $t \geqslant 0$. 
\end{definition}

\begin{definition}
We say that an equilibrium $x_0$ is uniformly Lyapunov stable if for any open set $U \subset M$ such that $x_0 \in U$, there exists an open set $V \subset M$ such that
$$
\Phi_{t_0, t}(x) \in U
$$
for any $x \in V$ and all $t_0 \geqslant 0$ and $t \geqslant 0$.
\end{definition}

The uniformity in the above definition is uniformity in time: the open neighborhood $V$ does not depend on $t_0$. In particular, if we consider a continuous semi-flow on $M$ (a continuous semi-process on $M \times \mathbb{R}$ without any dependence on $t_0$), then any Lyapunov stable equilibrium, defined in the usual way, is uniformly stable.
\begin{definition}
We say that $x_0$ is globally attractive if $\Phi_{0, t}(x) \to x_0$ as $t \to \infty$ for any $x \in M$. 
\end{definition}

\begin{definition}
Let $S \subset M \times \mathbb{R}$. Define a subset of $M \times \mathbb{R}$ by
$$
S_\tau = \{ (x, t) \in S \colon t = \tau \}.
$$
\end{definition}
In other words, $S_\tau$ is the section of $S$ by the plane $t = \tau$. For the sets of strict egress points, we will use the notation $W^{--}_0 = (W^{--})_0$.

We will now prove the main result. The proof is based on the idea of the Wa{\.z}ewski topological method.

\begin{theorem}
Let $W \subset M \times \mathbb{R}$ be an open set, $\Phi$ be a continuous semi-process on $M$ and $W^- = W^{--}$ (w.r.t. $\Phi$) and assume $W_0^{--} \ne \varnothing$. Let $x_0$ be a uniformly stable equilibrium, $U \subset M$ be an open subset, and $x_0 \in U$ be such that $\bar U \times \mathbb{R} \subset W$. Suppose that $W_0^{--}$ can be connected with the equilibrium by a continuous path $\Gamma \colon [0,1] \to M \times \mathbb{R}$ such that $\Gamma(s) \in W_0$ for $s \in (0,1)$ and $\Gamma(0) = (x_0, 0)$, $\Gamma(1) \in W_0^{--}$. Then $x_0$ cannot be globally attractive for $\Phi$.
\end{theorem}
\begin{proof}

We will prove by contradiction that there is a point $(x,0) \in \Gamma$ such that $\sigma(x) = \infty$ and $\Phi_{0, t}(x) \not\to x_0$ as $t \to \infty$.

For any point $(x,0) \in \Gamma$ we have two options: either the corresponding trajectory leaves $W$ ($\sigma(x) < \infty$), or the trajectory always remains inside $W$ ($\sigma(x) = \infty$). Note that neither of these sets are empty since $\Gamma(0) = (x_0, 0)$ corresponds to the equilibrium and $\Gamma(1) \in W_0^{--}$ and, therefore, the solution starting at this end of $\Gamma$, leaves $W$.

Assume that for all points $(x,0) \in \Gamma$ satisfying $\sigma(x) = \infty$ we have $\Phi_{0, t}(x) \to x_0$ as $t \to \infty$.

Now consider the following map $\Omega$ from $\Gamma$ to its boundary points $\Gamma(0)$ and $\Gamma(1)$:
\begin{align}
    \Omega(x, 0) = 
    \begin{cases}
    & \Gamma(0), \quad\mbox{if}\quad \sigma(x) = \infty,\\
    & \Gamma(1), \quad\mbox{if}\quad \sigma(x) < \infty. 
    \end{cases}
\end{align}

Now we will prove that $\Omega$ is continuous provided our assumption on the attractiveness holds.

Since the equilibrium $x_0$ is uniformly stable, there is an open set $V \subset M$, $x_0 \in V$ such that for any $t_0 \geqslant 0$ and any $x \in V$ for all $t > 0$ we have $\Phi_{t_0,t}(x) \in U$.

\begin{figure}[h!]
  \centering
    \includegraphics[width=0.45\textwidth]{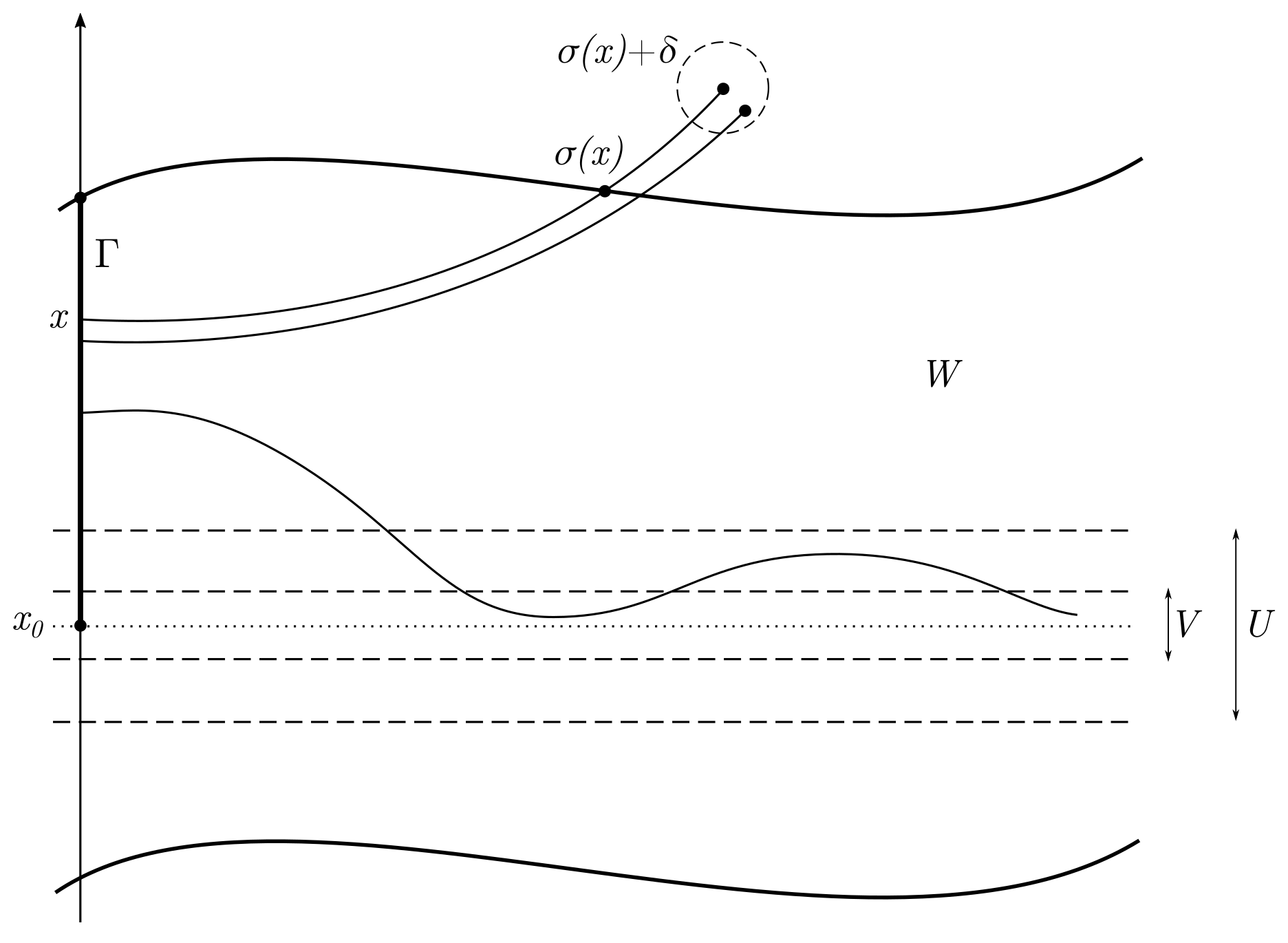}
\caption{A schematic representation of the map between $\Gamma$ and its boundary points.}
\label{fig4}
\end{figure}

If $\Omega(x,0) \mapsto \Gamma(0)$, then $\Omega(y,0) \mapsto \Gamma(0)$ provided $y$ is close to $x$: for some $\tau$ we have $\Phi_{0,\tau}(x) \in V$. Hence, $\Phi_{0,\tau}(y) \in V$ and $\Phi_{0, t}(y) \in U$ for all $t \geqslant \tau$ (Fig. \ref{fig4}). Therefore, the corresponding trajectory never leaves $W$ and $\sigma(y) = \infty$.

We will now prove that $\Omega$ is continuous at all points $(x,0)$ that are mapped to $\Gamma(1)$, i.e. $\sigma(x) < \infty$. Let $y$ be close to $x$. The point $(\Phi_{0,\sigma(x)}(x),\sigma(x))$ is a point of strict egress, therefore $(\Phi_{0,\sigma(x) + \delta}(x),\sigma(x) + \delta) \notin W\cup\partial W$ for some $\delta > 0$ (Fig. \ref{fig4}). Since the semi-process is continuous, we can conclude that $(\Phi_{0,\sigma(y) + \delta}(y),\sigma(y) + \delta)$ belongs to a small neighborhood of $(\Phi_{0,\sigma(x) + \delta}(x),\sigma(x) + \delta)$. In particular, $(\Phi_{0,\sigma(y) + \delta}(y),\sigma(y) + \delta) \notin W\cup\partial W$ and $\sigma(y)<\infty$.

Therefore, we have constructed a continuous map between a line segment and its boundary. From the contradiction we obtain that our assumption cannot be true and there exists a solution starting at $\Gamma$ such that this solution never leaves $W$ and does not tend asymptotically to $x_0$. In particular, $x_0$ cannot be globally attractive.

\end{proof} 

One can compare the above statement with system \eqref{exam3}. There is a stable equilibrium inside the small circle and we have shown that there is a solution that both does not leave the bigger circle and does not intersect the smaller one (Fig. 3). Note that the problem is not in the trajectories that leave our region and cannot return without breaking the continuity of the flow, but there always exists a solution that does not intersect our circles and is separated from the equilibrium. In other words, even if we make our system discontinuous and instantly carry over to the equilibrium all the solutions leaving the bigger circle, we still cannot obtain  global stabilization.

From the theorem we have that for each curve $\Gamma$ we have at least one point $x$ such that $\Phi_{0, t}(x) \not\to x_0$. Therefore, if we can find an $n$-parameter family of disjoint curves $\Gamma$, then we obtain an $n$-parameter family of the corresponding solutions.

Results similar to Theorem 2.13 can be proved for the case when our system has a stable invariant manifold. For instance, the following generalization can be considered.

Let $M = S \times N$ where $S$ and $N$ are smooth manifolds. If $x_0 \in S$, we say that $\{x_0\} \times N$ is an invariant manifold for $\Phi$ if 
for any $(x_0,y) \in \{x_0\} \times N$, $t_0 \in \mathbb{R}$ and $t \geqslant 0$ we have
$$
\Phi_{t_0, t}(x_0,y) \in \{x_0\} \times N.
$$
We say that an invariant manifold $\{x_0\} \times N$ is uniformly Lyapunov stable if for any open set $U \subset S$ such that $x_0 \in U$, there exists an open set $V \subset S$, $x_0 \in V$ such that
$$
\Phi_{t_0, t}(x,y) \in U \times N
$$
for any $x \in V$, $y \in N$ and all $t_0$ and $t \geqslant 0$.

We say that an invariant manifold $\{ x_0 \} \times N$ is globally attractive if for any $(x,y) \in S \times N$ we have $\Phi_{0, t} (x,y) \to \{ x_0 \} \times N$ as $t \to \infty$. We  will denote the canonical projection onto the manifold $M$ by $\pi_M \colon M \times \mathbb{R} \to M$.

The proof of the following result is the same as in Theorem 2.13. 

\begin{theorem}
Let $M = S \times N$, where $S$ and $N$ are smooth manifolds, and let $\Phi$ be a continuous semi-process on $M$. Let $W \subset M \times \mathbb{R}$ be such that  $W^- = W^{--}$ (w.r.t. $\Phi$) and assume $W_0^{--} \ne \varnothing$. Let $\{ x_0 \} \times N \subset M$ be a uniformly stable invariant manifold, $U \subset S$ be an open subset, $x_0 \in U$ and $\bar U \times N \times \mathbb{R} \subset W$. Suppose that $W_0^{--}$ can be connected with the invariant manifold by a continuous path $\Gamma \colon [0,1] \to M \times \mathbb{R}$ such that $\Gamma(s) \in W_0$ for $s \in (0,1)$ and $\pi_M(\Gamma(0)) \in \{ x_0 \} \times N$, $\Gamma(1) \in \partial W \cap W_0^{--}$. Then $\{ x_0 \} \times N $ cannot be globally attractive for $\Phi$.
\end{theorem}

\subsection{Remarks on the condition $W^- = W^{--}$}

It is also worth mentioning that the verification of the fact that $W^- = W^{--}$ can be simplified when our semi-process is defined by an ODE in a neighborhood of $\partial W$. In this case we can use definitions similar to the ones given above for flows.

To be more precise, suppose given a smooth manifold $M$ and a continuous semi-process $\Phi$ on it. Furthermore, suppose given an open set $W \subset M \times \mathbb{R}$. Suppose that $\Phi$ can be defined by a smooth ODE in a neighborhood of $\partial W$, i.e. there exists an open set $O(\partial W)$ such that $\partial W \subset O(\partial W)$ and there exists an ordinary differential equation
\begin{align}
\label{eq555}
    \dot x = v(x,t),
\end{align}
where $v(x,t)$ is a smooth function on $O(\partial W)$. Moreover, for any $(x_0,t_0) \in O(\partial W)$, $\Phi_{t_0,t}(x_0)$ is a differentiable function that satisfies \eqref{eq555}. In particular, the semi-process $\Phi$ can be considered in reverse time in $O(\partial W)$.

\begin{figure}[h!]
  \centering
    \includegraphics[width=0.45\textwidth]{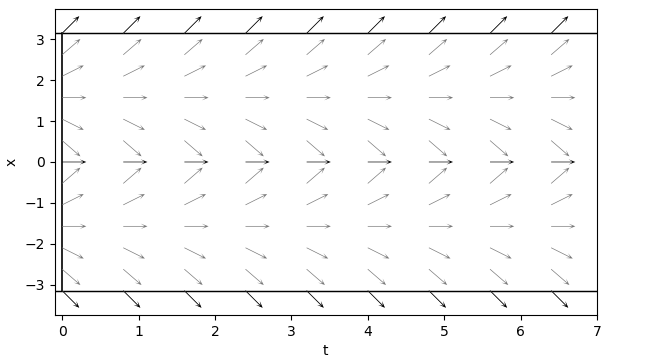}
\caption{A semi-process (semi-flow) defined by an ODE with a discontinuous right hand side.}
\label{fig5}
\end{figure}

\begin{definition}
We say that $(x_0,t_0) \in \partial W$ is an egress point for $W$ w.r.t. \eqref{eq555} if there exists an $\varepsilon > 0$ such that $(\Phi_{t_0,t}(x_0), t_0+t) \in W$ for all $t \in (-\varepsilon,0)$.
\end{definition}
\begin{definition}
We say that an egress point $(x_0,t_0) \in \partial W$ is a strict egress point for $W$ w.r.t. \eqref{eq555} if there exists an $\varepsilon > 0$ such that $(\Phi_{t_0,t}(x_0), t_0+t) \notin W \cup \partial W$ for all $t \in (0,\varepsilon)$.
\end{definition}

Therefore, if $W^{-} = W^{--}$ in the sense of Definitions 2.15 and 2.16, then $W^{-} = W^{--}$ in the sense of Definitions 2.6 and 2.7. Note that in the latter case, the set of egress points is a subset of $W^-$ in the sense of Definition 2.15.

As an illustration of this approach, consider the system (Fig. 5)
$$
    \dot x = 
    \begin{cases}
    -\cos(x) &\mbox{if } x > 0, \\
0 & \mbox{if } x = 0,\\
\cos(x)  &\mbox{if } x < 0.
    \end{cases}
$$
This system defines a semi-process (semi-flow): for $x_0 = 0$ we cannot consider the corresponding solution $x = 0$ for $t < t_0$. There are infinitely many possibilities for a continuous continuation. However, if we consider the set $W = \{ x \colon -\pi < x <\pi \}$, we see that in a neighborhood of the boundary our semi-flow is defined by the solutions of two ODEs: $\dot x = -\cos(x)$ and $\dot x = \cos(x)$. In the sense of Definitions 2.15 and 2.16, we have $W^- = W^{--}$. Therefore, any solution starting from the interval $t = 0$, $-\pi \leqslant x \leqslant \pi$ that can reach the boundary, locally leaves $W \cup \partial W$ and $W^- = W^{--}$ in the sense of Definitions 2.6 and 2.7. Moreover, $x = 0$ is a uniformly stable equilibrium and we can apply Theorem 2.13. There exists a solution ($x = \pi/2$ or $x = -\pi/2$) that never leaves $W$ and does not tend asymptotically to an equilibrium.

\section{Examples}

\subsection{The inverted pendulum}
We will begin with the equation
\begin{align}
\label{eq5}
    \ddot \varphi = u(\varphi, \dot \varphi) \sin \varphi - \cos \varphi + v(\varphi, \dot \varphi).
\end{align}

Equation (\ref{eq5}) describes the motion of a controlled inverted pendulum in a gravitational field. The feedback control is given by $u$ and $v$. The function $u$ defines the horizontal acceleration of the pivot point and $v$ is a control torque.

First, we assume that this equation defines a continuous flow on $\mathbb{R}^2$. Here $u, v \colon \mathbb{R}^2 \to \mathbb{R}$ are smooth. We also assume that $|v(0,0)| < 1$ and $|v(\pi,0)| < 1$, with $\varphi = \pi/2$ being a Lyapunov stable equilibrium (this equilibrium can be made stable by choosing an appropriate $u$ and $v$). Then this equilibrium cannot be globally attractive.

For this system, $W$ has the following form
$$
W = \{ \varphi, \dot\varphi, t \colon 0 < \varphi < \pi \}.
$$
We will show that $W^- = W^{--}$. This follows from the Taylor expansion for $\varphi(t)$. Indeed, let $\varphi_0 = 0$ and $\dot \varphi_0 < 0$. Then
$$
\varphi(t) = \varphi_0 + (t - t_0) \dot \varphi_0 + o(|t - t_0|).
$$
Therefore, we can conclude that $(\varphi_0, \dot \varphi_0, t_0) \in W^-$ (i.e., $\varphi(t) > 0$ for $t < t_0$ provided $|t - t_0|$ is small). In accordance with Definition 2.15, $(\varphi_0, \dot\varphi_0,t_0)$ is an egress point. Similarly, $(\varphi_0, \dot \varphi_0, t_0) \in W^{--}$ (i.e., $\varphi(t) < 0$ for $t > t_0$ provided $|t - t_0|$ is small). When $\varphi_0 = 0$ and $\dot\varphi_0 = 0$, we have 
$$
\varphi(t) = \frac{1}{2}(t - t_0)^2 \ddot \varphi (t_0) + o(|t - t_0|^2) = \frac{1}{2}(t - t_0)^2 (v(0,0) - 1) + o(|t - t_0|^2).
$$
Therefore, $(\varphi_0, \dot \varphi_0, t_0) \not\in W^-$. Finally, one can show that
$$
W^- = W^{--} = \{ \varphi, \dot\varphi, t \colon \varphi = 0, \dot\varphi < 0 \mbox{ or } \varphi = \pi, \dot\varphi > 0\}.
$$

Similar considerations can be found, for instance, in \cite{polekhin2014periodic,polekhin2014examples}, where they were used to prove the existence of non-falling and periodic solutions for the pendulum with a moving pivot point.

There is a solution of (\ref{eq5}) whose trajectory always remains in $W$ and does not tend to the equilibrium. Therefore, $(\pi/2,0)$ cannot be a globally attractive uniformly stable equilibrium in the system where the pendulum moves along the plane of support (the horizontal line) and the rod can hit this plane. Taking into account the fact that the trajectory of the solution always remain inside the set $0 < \varphi < \pi$, i.e., along this solution the rod of the pendulum never becomes horizontal, we can conclude that the existence of this solution does not depend on the model of impact between the rod and the horizontal plane.

Note that the same result holds for the system
$$
    \ddot \varphi = u(\varphi, \dot \varphi, t) \sin \varphi - \cos \varphi + v(\varphi, \dot \varphi, t) + f(t)\sin\varphi,
$$
where $f(t)$ is an external horizontal force acting on the pendulum (a smooth function). For this system we also have $W^- = W^{--}$ for the same $W$.

A similar result can be obtained for the following system describing the motion of an inverted pendulum on a cart. The details can be found in \cite{polekhin2018topological}.

\subsection{The Furuta pendulum}
We will now consider the Furuta pendulum, a well known control system introduced in \cite{furuta1991swing} and thoroughly studied by many authors (for instance, \cite{furuta1992swing,shiriaev2007virtual,cazzolato2011dynamics,ramirez2014linear,nair2002normal,la2009new,aracil1998global,ibanez2007stabilization,shiriaev2001stabilization}).
 The system consists of a pendulum and a rotating base on which the pendulum is mounted. The pendulum is controlled by a torque applied to the base (Fig. 6). The governing equations have the form

\begin{align}
\begin{split}
\label{eqq8}
    &[I + m(L + l^2 \sin^2 \varphi)] \ddot \psi + mlL \cos \varphi \ddot \varphi + \frac{1}{2}ml^2 \dot \varphi \sin 2\varphi \dot \psi +\\
    &\left[-mlL\sin\varphi + \frac{1}{2}ml^2 \sin 2\varphi \dot \psi \right] \dot\varphi = u(\varphi, \dot \varphi, \psi, \dot\psi),\\
    & mlL \cos\varphi \ddot\psi + ml^2  \ddot\varphi -\frac{1}{2} ml^2 \sin 2\varphi \dot\psi \dot\varphi - mgl\sin\varphi = 0.
\end{split}
\end{align}

\begin{figure}[h!]
  \centering
    \def\svgwidth{250 pt}
\begingroup%
  \makeatletter%
  \providecommand\color[2][]{%
    \errmessage{(Inkscape) Color is used for the text in Inkscape, but the package 'color.sty' is not loaded}%
    \renewcommand\color[2][]{}%
  }%
  \providecommand\transparent[1]{%
    \errmessage{(Inkscape) Transparency is used (non-zero) for the text in Inkscape, but the package 'transparent.sty' is not loaded}%
    \renewcommand\transparent[1]{}%
  }%
  \providecommand\rotatebox[2]{#2}%
  \newcommand*\fsize{\dimexpr\f@size pt\relax}%
  \newcommand*\lineheight[1]{\fontsize{\fsize}{#1\fsize}\selectfont}%
  \ifx\svgwidth\undefined%
    \setlength{\unitlength}{579.8539866bp}%
    \ifx\svgscale\undefined%
      \relax%
    \else%
      \setlength{\unitlength}{\unitlength * \real{\svgscale}}%
    \fi%
  \else%
    \setlength{\unitlength}{\svgwidth}%
  \fi%
  \global\let\svgwidth\undefined%
  \global\let\svgscale\undefined%
  \makeatother%
  \begin{picture}(1,0.4798567)%
    \lineheight{1}%
    \setlength\tabcolsep{0pt}%
    \put(0,0){\includegraphics[width=\unitlength,page=1]{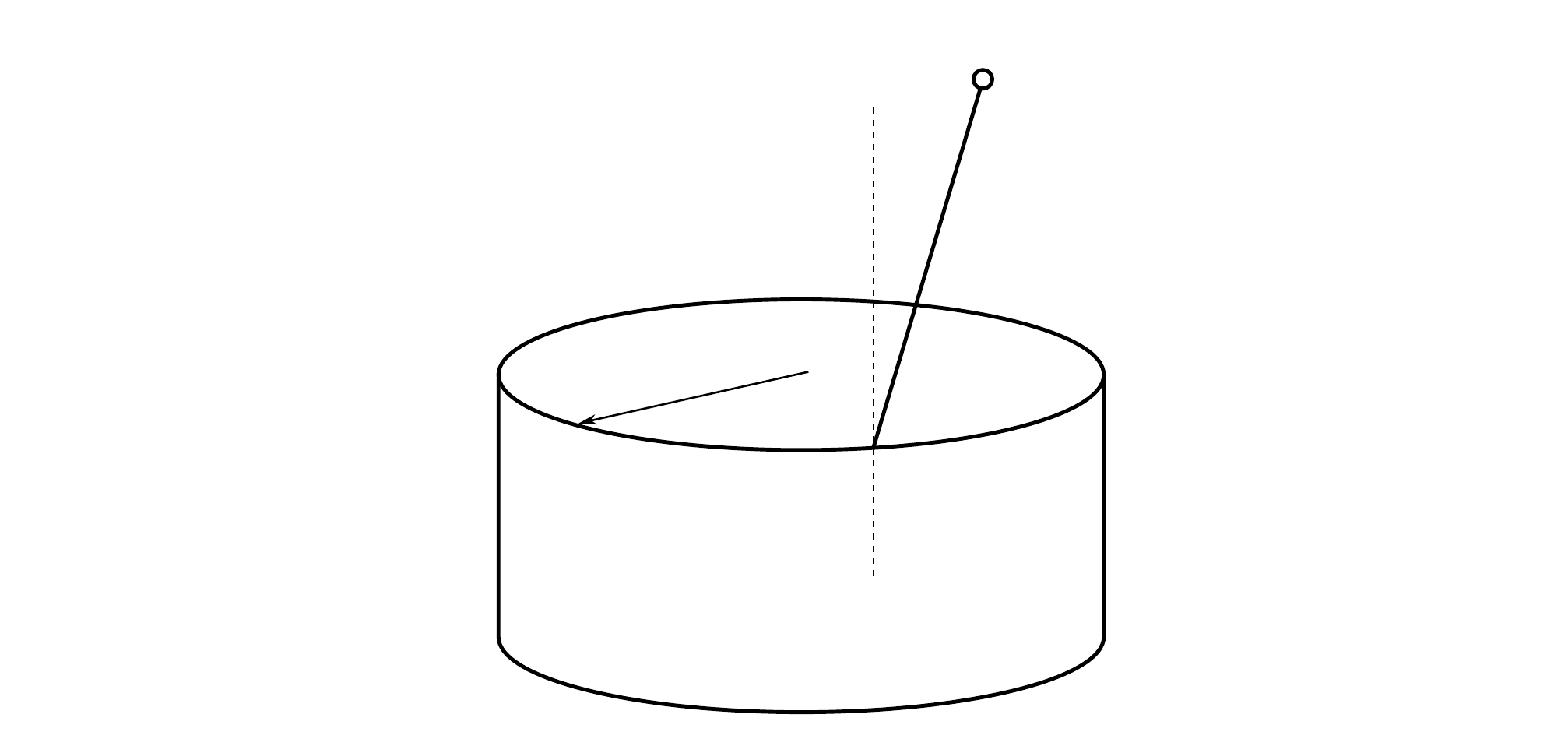}}%
    \put(0.41518798,0.23623285){\color[rgb]{0,0,0}\makebox(0,0)[lt]{\lineheight{1.25}\smash{\begin{tabular}[t]{l}$L$\end{tabular}}}}%
    \put(0.64504885,0.42581258){\color[rgb]{0,0,0}\makebox(0,0)[lt]{\lineheight{1.25}\smash{\begin{tabular}[t]{l}$m$\end{tabular}}}}%
    \put(0,0){\includegraphics[width=\unitlength,page=2]{scl-fig3.pdf}}%
    \put(0.64776032,0.07332116){\color[rgb]{0,0,0}\makebox(0,0)[lt]{\lineheight{1.25}\smash{\begin{tabular}[t]{l}$M$\end{tabular}}}}%
    \put(0.62152203,0.36249494){\color[rgb]{0,0,0}\makebox(0,0)[lt]{\lineheight{1.25}\smash{\begin{tabular}[t]{l}$l$\end{tabular}}}}%
    \put(0,0){\includegraphics[width=\unitlength,page=3]{scl-fig3.pdf}}%
    \put(0.27960029,0.37727699){\color[rgb]{0,0,0}\makebox(0,0)[lt]{\lineheight{1.25}\smash{\begin{tabular}[t]{l}$g$\end{tabular}}}}%
    \put(0.57043706,0.35436483){\color[rgb]{0,0,0}\makebox(0,0)[lt]{\lineheight{1.25}\smash{\begin{tabular}[t]{l}$\varphi$\end{tabular}}}}%
    \put(0,0){\includegraphics[width=\unitlength,page=4]{scl-fig3.pdf}}%
  \end{picture}%
\endgroup%

\caption{Furuta pendulum: An inverted pendulum mounted on a rotating base.}
\label{fig:fig}
\end{figure}

For the sake of brevity, we assume that the system is moving without any friction. Here $I$ is the inertia of the rotating base, $L$ is the radius of the base, $l$ is the length of the pendulum, $m$ and $M$ are the masses of the pendulum and the base, respectively. $\varphi$ is the angle between the upward vertical direction and the rod of the pendulum, $\psi$ is the angle of rotation of the base. We do not assume that $u$ is periodic in $\varphi$ or $\psi$, i.e. $u \colon \mathbb{R}^4 \to \mathbb{R}$. In particular, we cannot apply here the result from \cite{bhat2000topological} since we do not have any rotational degrees of freedom.

Now assume that the solutions define a semi-flow on $\mathbb{R}^4$ and $u$ is a smooth function in a neighborhood of points where $\varphi = \pi/2$ or $\varphi = -\pi/2$. Then the point $\varphi = 0$, $\psi = \psi_0$ cannot be a globally attractive Lyapunov stable equilibrium. Indeed, let us consider the set
$$
W = \{ \varphi, \dot \varphi, \psi, \dot\psi, t \colon -\pi/2 < \varphi < \pi/2 \}
$$
The right hand side of the system is smooth in a neighborhood of the boundary $\partial W$. We can consider Taylor expansions similar to those presented above for the case of an inverted pendulum. From the second equation of system \eqref{eqq8} we have $\ddot\varphi > 0$ when $\dot\varphi = 0$ and $\varphi = \pi/2$ and $\ddot\varphi < 0$ when $\dot\varphi = 0$ and $\varphi = -\pi/2$. Therefore $W^- = W^{--}$ and Theorem 2.13 can be applied.

\subsection{The wheeled pendulum}

Another example related to the dynamics of pendulum-like systems is the wheeled pendulum. For instance, various results for one- and two-wheeled pendulums can be found in \cite{formalsky2014motion,pathak2005velocity,li2010robust}. Note that the wheeled pendulum can be considered as a model for Segway, a self-balancing personal transporter \cite{lee2008control,do2010motion}. 

The wheeled pendulum is a mathematical pendulum with its pivot point attached to a disk rolling without slipping on a horizontal line (Fig. 7). We consider this system as a controlled system and assume that there is a control torque $u$ applied to the pivot.

\begin{figure}[h!]
  \centering
    \def\svgwidth{250 pt}
\begingroup%
  \makeatletter%
  \providecommand\color[2][]{%
    \errmessage{(Inkscape) Color is used for the text in Inkscape, but the package 'color.sty' is not loaded}%
    \renewcommand\color[2][]{}%
  }%
  \providecommand\transparent[1]{%
    \errmessage{(Inkscape) Transparency is used (non-zero) for the text in Inkscape, but the package 'transparent.sty' is not loaded}%
    \renewcommand\transparent[1]{}%
  }%
  \providecommand\rotatebox[2]{#2}%
  \newcommand*\fsize{\dimexpr\f@size pt\relax}%
  \newcommand*\lineheight[1]{\fontsize{\fsize}{#1\fsize}\selectfont}%
  \ifx\svgwidth\undefined%
    \setlength{\unitlength}{579.8539866bp}%
    \ifx\svgscale\undefined%
      \relax%
    \else%
      \setlength{\unitlength}{\unitlength * \real{\svgscale}}%
    \fi%
  \else%
    \setlength{\unitlength}{\svgwidth}%
  \fi%
  \global\let\svgwidth\undefined%
  \global\let\svgscale\undefined%
  \makeatother%
  \begin{picture}(1,0.4798567)%
    \lineheight{1}%
    \setlength\tabcolsep{0pt}%
    \put(0,0){\includegraphics[width=\unitlength,page=1]{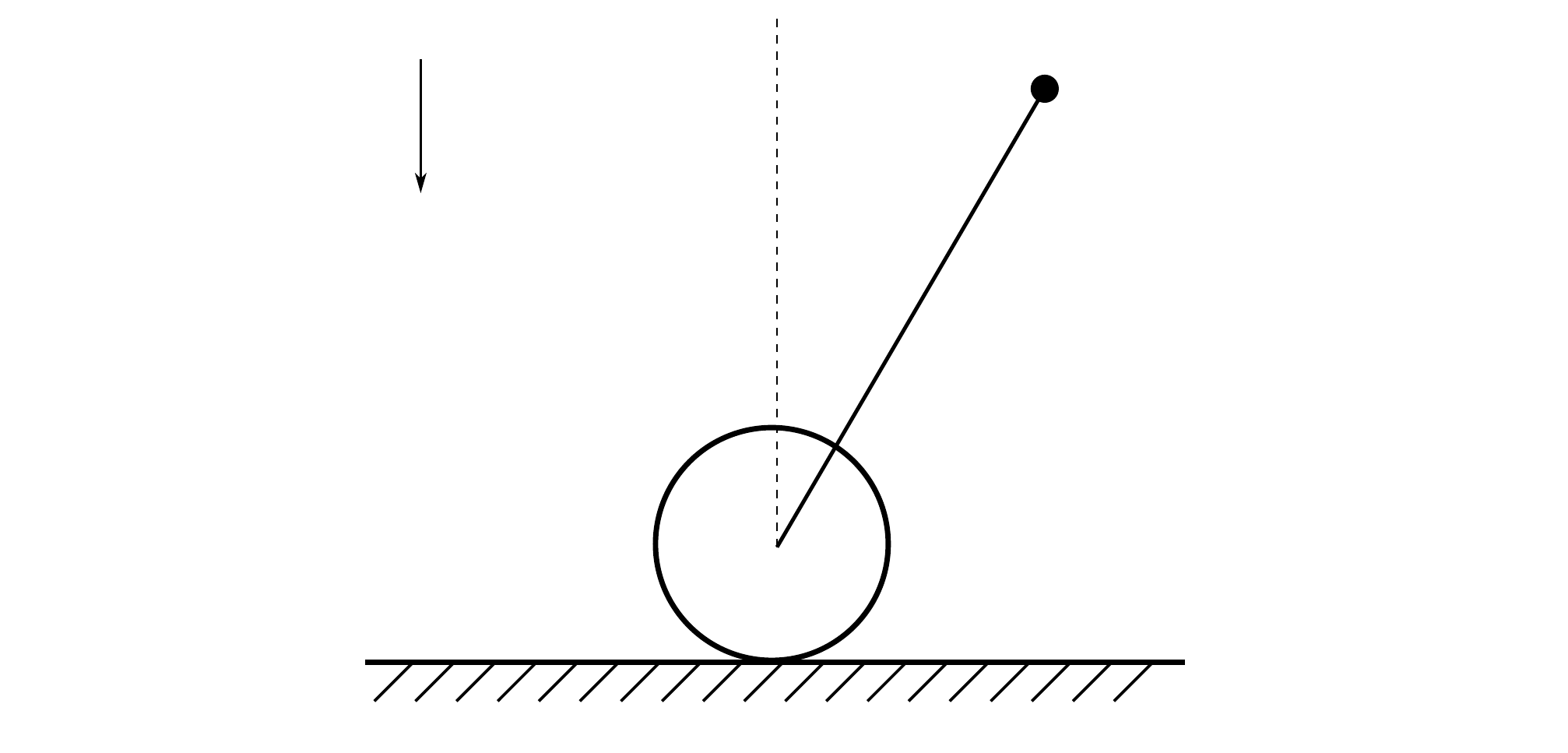}}%
    \put(0.27995744,0.39238441){\color[rgb]{0,0,0}\makebox(0,0)[lt]{\lineheight{1.25}\smash{\begin{tabular}[t]{l}$g$\end{tabular}}}}%
    \put(0.58004219,0.12540873){\color[rgb]{0,0,0}\makebox(0,0)[lt]{\lineheight{1.25}\smash{\begin{tabular}[t]{l}$M$\end{tabular}}}}%
    \put(0.68751552,0.41563968){\color[rgb]{0,0,0}\makebox(0,0)[lt]{\lineheight{1.25}\smash{\begin{tabular}[t]{l}$m$\end{tabular}}}}%
    \put(0,0){\includegraphics[width=\unitlength,page=2]{scl-fig2.pdf}}%
    \put(0.59434793,0.26221485){\color[rgb]{0,0,0}\makebox(0,0)[lt]{\lineheight{1.25}\smash{\begin{tabular}[t]{l}$l$\end{tabular}}}}%
    \put(0.52810835,0.28722992){\color[rgb]{0,0,0}\makebox(0,0)[lt]{\lineheight{1.25}\smash{\begin{tabular}[t]{l}$\varphi$\end{tabular}}}}%
    \put(0,0){\includegraphics[width=\unitlength,page=3]{scl-fig2.pdf}}%
    \put(0.44148404,0.13493681){\color[rgb]{0,0,0}\makebox(0,0)[lt]{\lineheight{1.25}\smash{\begin{tabular}[t]{l}$r$\end{tabular}}}}%
    \put(0,0){\includegraphics[width=\unitlength,page=4]{scl-fig2.pdf}}%
  \end{picture}%
\endgroup%

\caption{An inverted pendulum on a wheel.}
\label{fig:fig}
\end{figure}

By $m$ and $M$ we denote the masses of the pendulum and the disk, respectively. Let $l$ be the length of the pendulum and $r$ is the radius of the disk. It can be shown that the equation for the angle between the rod and the vertical direction can be considered independently and has the form
\begin{align*}
    \begin{split}
        &(a_{11}a_{22} - a_{12}^2\cos\varphi)\ddot\varphi + a_{12}^2 \dot\varphi^2 \sin\varphi\cos\varphi -\\
        &a_{11}mgl\sin\varphi = (a_{11}+a_{12}\cos\varphi) u(\varphi, \dot\varphi,t).
    \end{split}
\end{align*}
Here $a_{11} = (2M + m)r^2$, $a_{12} = mrl$, $a_{22} = ml^2$, while $u \colon \mathbb{R}^3 \to \mathbb{R}$ is smooth in a vicinity of the planes $\varphi = \pi/2$ and $\varphi = -\pi/2$. The function $u$ defines the control torque applied to the rod. We assume that solutions of the above equation define a semi-process on $\mathbb{R}^2$. Then the vertical upward position cannot be a globally attractive uniformly stable equilibrium provided $|u| < mgl$ holds for all $t$, and for $(\varphi, \dot\varphi) = (\pi/2,0)$ and $(\varphi, \dot\varphi) = (-\pi/2,0)$. Again, here we have $\ddot\varphi > 0$ when $\dot\varphi = 0$ and $\varphi = \pi/2$ and $\ddot\varphi < 0$ when $\dot\varphi = 0$ and $\varphi = -\pi/2$ and for $W = \{ \varphi, \dot\varphi, t \colon -\pi/2 < \varphi < \pi/2 \}$ we obtain $W^- = W^{--}$.

\section{Conclusion}
 We have presented an approach that can be applied to many real-life systems when one wants to prove that the system cannot be globally stabilized by means of a feedback control. There are two key requirements that need to be met to apply the method. First, the right hand sides of the equations should be relatively regular functions. At least, they should define a semi-flow (or a semi-process, if the system is non-autonomous). For instance, in the above examples we only assumed that the functions that define our feedback control are smooth only in a neighborhood of the boundary of $W$. Second, we need a subset of the extended phase space such that, broadly speaking, all trajectories of our control system are either transverse to the boundary of this set, or externally tangent to it and at least some trajectories leave the region. Here it is worth mentioning that in this case the reason the system cannot be globally stabilized is not because of the solutions that leave our set and cannot return to the equilibrium, but because there exists a solution (often it is a family of solutions) that neither leaves our set nor asymptotically tends to the equilibrium.
 
As it is illustrated by the above examples, for various pendulum-like systems it is often possible to find the required set and to prove the impossibility of global stabilization. Note that it is usually needed to stabilize the upward vertical position in such systems. For this case it is possible to apply our results. For the downward vertical position the situation is different. For instance, for the mathematical pendulum \eqref{eq5}, the position $\varphi = -\pi/2$ can be globally asymptotically stable provided the rod of pendulum is constrained: $\varphi \in [-\pi,0]$. If we assume that the collisions between the constraint and the rod are elastic, then one can put here $\nu = -\dot\varphi$. One can add that it is also possible to prove the impossibility of global stabilization for the Lagrange top (3D pendulum) controlled by horizontal forces \cite{polekhin2018impossibility}.

In some sense the above results complement a theorem of Bhat and Bernstein \cite{bhat2000topological}. Our results are not so universal, but they can be applied to systems with non-compact configuration spaces or configuration spaces with boundaries and to non-autonomous systems.

\section*{Acknowledgement}
This work was performed at the Steklov International Mathematical Center and supported by the Ministry of Science and Higher Education of the Russian Federation (agreement no. 075-15-2019-1614).

\bibliographystyle{elsarticle-num}

\bibliography{sample}

\end{document}